# A Mathematical Theory for Studying and Controlling the Disinformation System Dynamics

Arindam Kumar Paul[1*] and M. Haider Ali Biswas[2]

[1]Mathematics Discipline, Science, Engineering and Technology School,
Khulna University, Khulna-9208, Bangladesh.
`1arindam017@gmail.com, 2mhabiswas@yahoo.com`

**Abstract**: This study explores the connection between disinformation, defined as deliberate spread of false information, and rate-induced tipping (R-tipping), a phenomenon where systems undergo sudden changes due to rapid shifts in external forces. While traditionally, tipping points were associated with exceeding critical thresholds, R-tipping highlights the influence of the rate of change, even without crossing specific levels. The study argues that disinformation campaigns, often organized and fast-paced, can trigger R-tipping events in public opinion and societal behavior. This can happen even if the disinformation itself doesn't reach a critical mass, making it challenging to predict and control. Here, by Transforming a population dynamics model into a network model, Investigating the interplay between the source of disinformation, the exposed population, and the medium of transmission under the influence of external sources, the study aims to provide valuable insights for predicting and controlling the spread of disinformation. This mathematical approach holds promise for developing effective countermeasures against this increasingly prevalent threat to public discourse and decision-making.

## 1  Introduction

Difference between Disinformation and Misinformation is in the intent, Disinformation always have clear defined intention but Misinformation may have intentions or may not. The purposefulness makes disinformation very complex, well-organized and chaotic. Reflexive Control theory initially was developed in Russia for use as a war strategy by Russia but applied to distract the European Union during the Ukraine invasion. Understanding the relationship between disinformation and public opinion is challenging with existing technology due to several reasons. The use of advanced technologies like artificial intelligence (AI) to spread disinformation poses additional challenges. For example, AI can be used to create realistic-looking fake accounts and content, making disinformation harder to detect[1], [2]. Often, the public's ability to discern false information online is limited, making them more susceptible to disinformation. Improving digital literacy is a long-term solution but poses its challenges. Harvard study rigorously stated the failure of AI in identifying and controlling Disinformation and recommended the necessity of inventing some out-of-the-box strategies[1], [3]. To do so, here the previously developed models to understand the disinformation and misinformation process by Paul & Biswas1 are used as the base dynamical system[4], [5]. In these previous two simultaneous studies on misinformation by Paul & Biswas, a new model was proposed first and developed considering the real cases and all the model was justified analytically and numerically. The significant findings of the study found a critical level for some transitions and successfully concluded by identifying the key factors of misinformation with and without intents. The study also proposes a rule-based solution for identifying the intention of the information spread[5].

## 2  Mathematical Model to Combat Disinformation

In this study, a previous model and analysis[4] have been extended for a better understanding and prediction. So, considering the Spreaders as the Medium or Target population by the source because they're very easy to influence, intended populations as the Sources who want to influence others by changing their beliefs with various strategies, Susceptible population as the initially connected population of a network, Exposed as the vulnerable to the Sources and the secondary targets of the sources[6]–[8]. Skeptics are being the Protesters who willingly think, unveil and share the truth and also creates a strong resistance powered by their higher level of authenticity. Another assumption is the rate of influences applied to the medium is a function of $\alpha$ and time t, the influences is turning Exposed to Medium transition rate and Susceptible to Medium transition rate dynamic. In this study, the precious model and analysis have been extended for a better understanding and prediction. So, extending the previous model can be written as:

$$\dot{x} = f(y, H(\alpha t)) \tag{1}$$

where, $\alpha \geq 0$ is a constant rate parameter. And consider when $t$ is the time scale of the system then $\tau = \alpha t$ as the time scale of the external input. Otherwise, instead of assuming $\tau = \alpha t$ and $H = f(\alpha, t)$, we can consider $H = f(\alpha, t, p)$, and the whole systems as

$$\dot{x} = f(y, H(\alpha, t, p)) \tag{2}$$

---
[1]The studies describe the patterns formation techniques of Misinformation spread using a system of nonlinear Differential Equations consist of 5 ODEs.

In other words, the rate parameter $\alpha$ quantifies the rate of change of external forcing with a given profile. By considering $\lambda(\epsilon, t)$ and $\Phi(\gamma, t)$ both are varying with the time for the influences of the time-varying rate of external influence and both are the functions of time $t$ [9], [10]. The whole system has been visualized in the Figure 1.

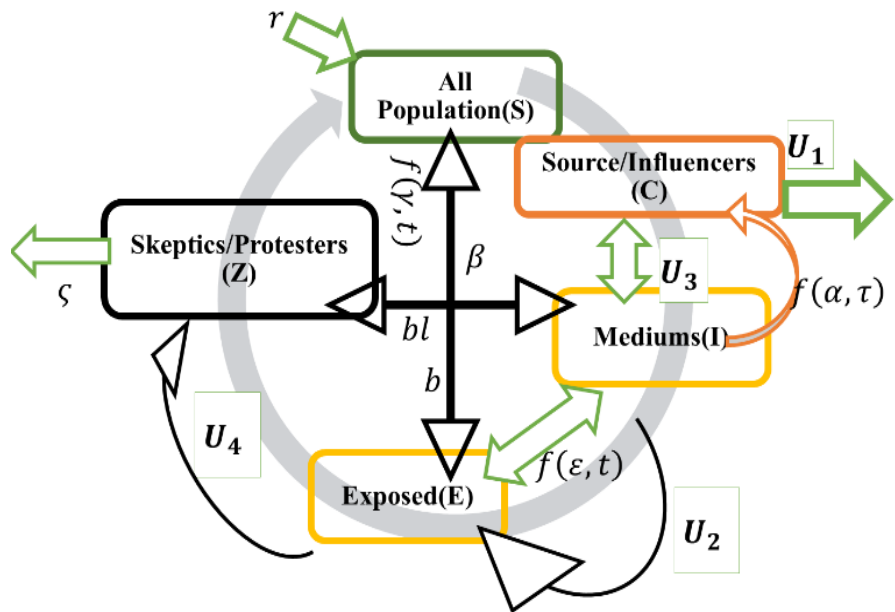

*Figure 1: A network diagram with 5 types of individuals representing the interaction and transition*

**Central Role of Mediums:** The model highlights the crucial role of Mediums (I) in mediating the spread of information or disinformation.

**Population Heterogeneity:** The model acknowledges that not all individuals respond identically to information or disinformation, with some becoming Exposed (E) and others becoming Skeptics/Protesters (Z).

The non-linear system of five ODEs are described in (3)-(7):

$$S' = r - f(\gamma, C, S) - \beta IS - bSZ \tag{3}$$

$$E' = \frac{-\varepsilon E}{\varphi Z} + (1-p)\beta IS + (1-l)bSZ \tag{4}$$

$$C' = f(\alpha, \tau)CI - \mu C + (1-\eta)f(\gamma, C, S) \tag{5}$$

$$I' = \frac{\varepsilon E}{\varphi Z} - f(\alpha,\tau)CI + \eta f(\gamma,C,S) + \beta pIS \qquad (6)$$

$$Z' = blZ - \xi Z \qquad (7)$$

Now after introducing the $f(\alpha,\tau)$ as the external forces responsible for the transmission of spreaders to intended for various influences, we'll examine the critical rate of change of $f(\alpha,\tau)$ in order to determine the optimal controlling strategies and the timeline of the application of control measures. The numerical simulation using all necessary parameters and initial conditions [2] have been visualized and described below:

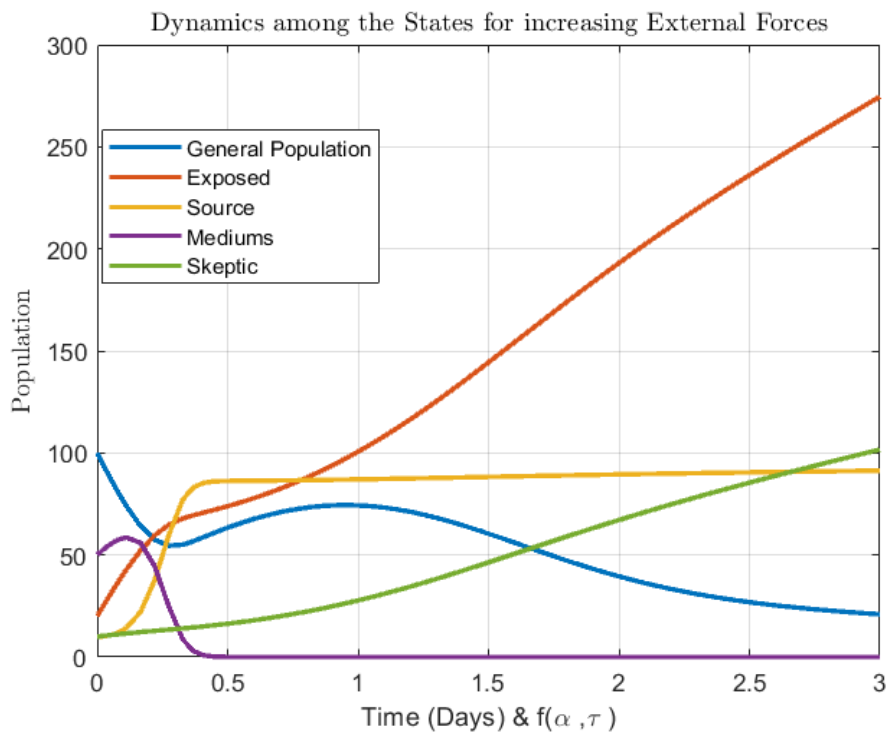

*Figure 2: The impact of external influence on the changes of system's components.*

From Figure 2, we find some significant insight into the system. The "Exposed" and "Intended" states are most affected by external forces, while "Spreaders" show moderate growth and "Skeptics" remain unaffected. Sharply increases with external forces, indicating many intend to spread disinformation.

---

[2] Description, Units and Numerical values or the parameters, initial conditions and others can be found at doi.org/10.13140/RG.2.2.29975.14243/1

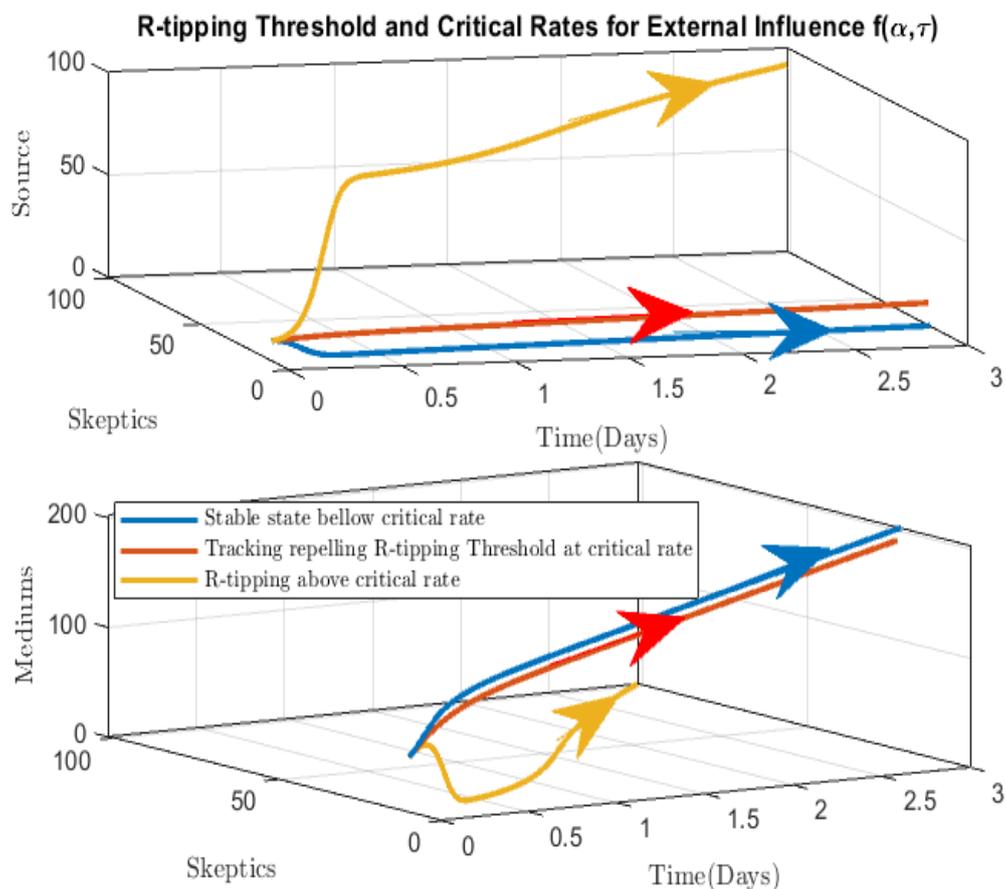

*Figure 3:,The interplay among the external influences f(α, τ)), Skeptics, Mediums and Source/Intended populations while Source tries to influence mediums to join them.*

In Figure 3, rapidly increasing external forces can induce complex and potentially irreversible changes in the Source and Medium population dynamics. Rapid acceleration of external forces Increase branching behavior of both due to the combined influences of external facts and skeptics, push system towards critical thresholds, Accelerate Medium conversion to Sources and Trigger sudden shifts in Medium population.

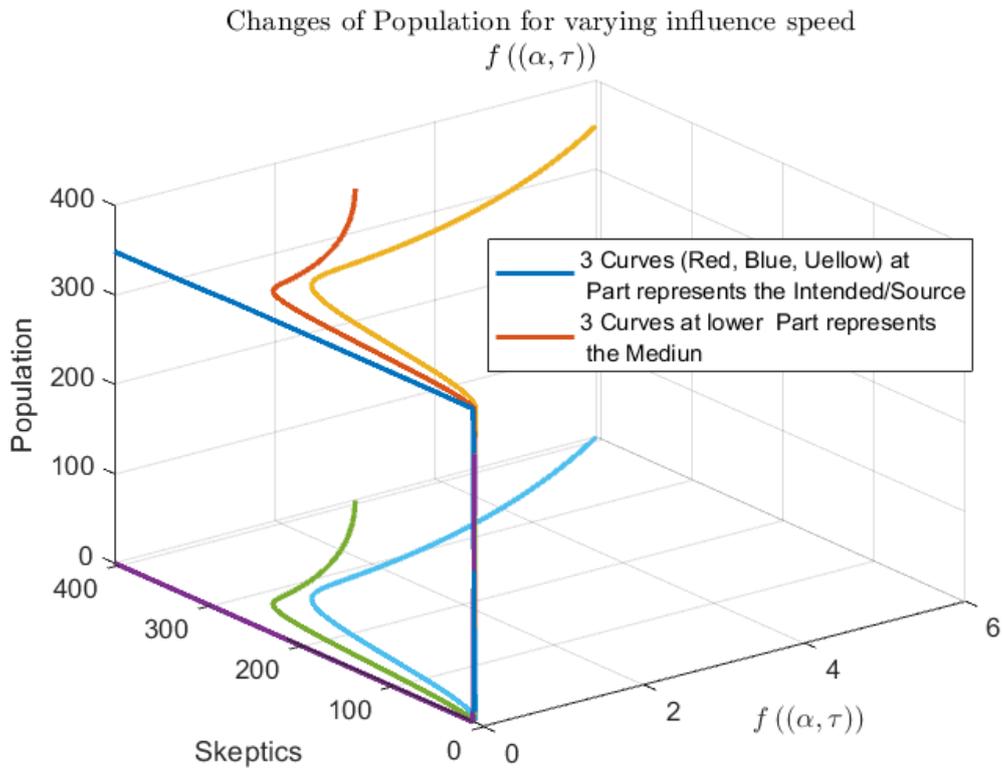

*Figure 4: the interplay between external influences (f(α, τ)) and the time-varying transmission rate (f(ε,t)) creates a complex and dynamic system with the potential for rapid changes and diverse outcomes. Formally it's well-known as Rate induced Tipping.*

Rapid acceleration of external forces can: Increase branching behavior of f(α, τ) This potentially leads to diverse outcomes depending on the chosen pathway within the system. While the fastest external influence represented by the red curve indicating that it continues increasing by a certain concerning rate even after the increment of skeptics' population. With the low rate of external influence represented by the yellow curve indicating that with the increase of skeptics' population Suddenly prevent the high rate of increment.

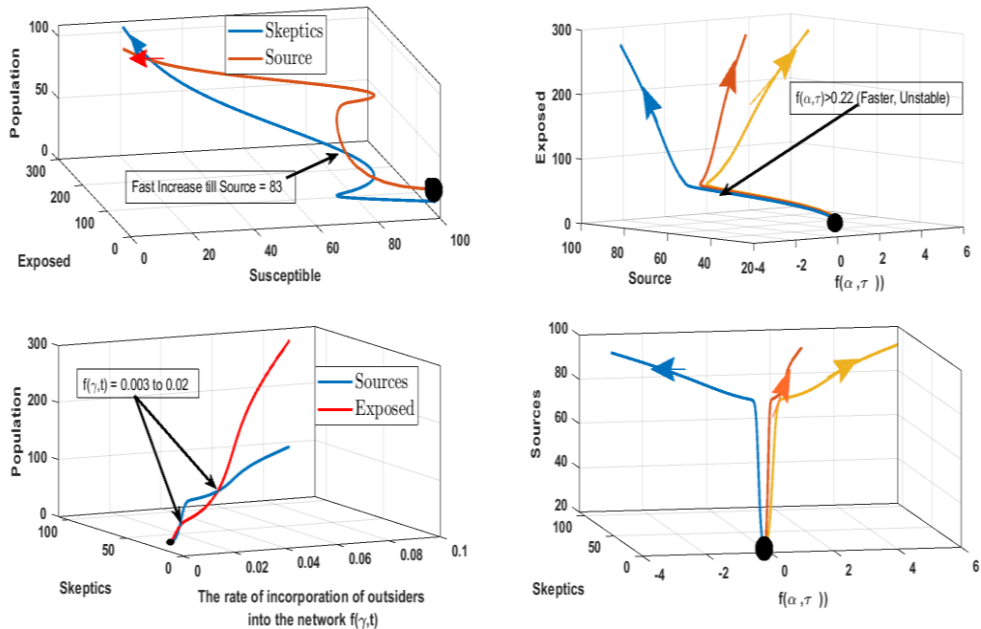

*Figure 5: The complex dynamical interplay among the external influences, skeptics, time-varying transmission rate $f(\gamma, t)$, and other population*

Figure 5 shows that the intended/source decreases with both external influences and skeptics increments, highlighting the need for interventions against both factors. Interaction among Source Population, Skeptics, and Transmission Rate $f(\gamma, t)$ indicates the properties of skeptics to freeze the impact of external forces on the system after a certain period of time, slow impact of influences on $f(\gamma, t)$. In the 1st Figure, time-varying transmission rate $f(\gamma, t)$ showcasing its changes and branching behavior under the influence of rapid acceleration of external forces (f(α, τ)). Fast Increase till Source = 83" near the Exposed highlights a potential tipping point where the rate of skepticism increases rapidly when the source reaches a critical value of 83. f(α, τ) > 0.22 suggests that rapid changes in influence (high α) can lead to unpredictable dynamics and potentially accelerate the spread of disinformation.

## 3    Application of Optimal Control in Policy Making

According the findings of simulations, we deployed four control measures as described in Figure 1 by $u_1, u_2, u_3, u_4$. The Objective functional is defined by

$$\min \int_0^T f(t, x(t), u(t))\, dt,$$

$$\begin{cases} f(t, x, u) = C + I + E - Z + u_1^2 + u_2^2 + u_3^2 + u_4^2 \\ x(t) = (S, E, C, I, Z) \\ S', E, C, I, Z' \in g(t, x, u) = x'(t) \end{cases} \qquad (8)$$

The existence and uniqueness of the control can be found and analytically represented by the 5.

$$u_1^* = \min\left\{\max\left(0, \frac{C\lambda_3}{2}\right), 1\right\}, u_2^* = \min\left\{\max\left(0, \frac{I\lambda_4}{2} - \frac{I\lambda_2}{2}\right), 1\right\} \qquad (9)$$

$$u_3^* = \min\left\{\max\left(0, \frac{CI\lambda_4}{2} - \frac{CI\lambda_3}{2}\right), 1\right\}, u_4^* = \max\left\{\min\left(0, \frac{SZb\lambda_5}{2} - \frac{SZb\lambda_2}{2}\right), 1\right\}$$

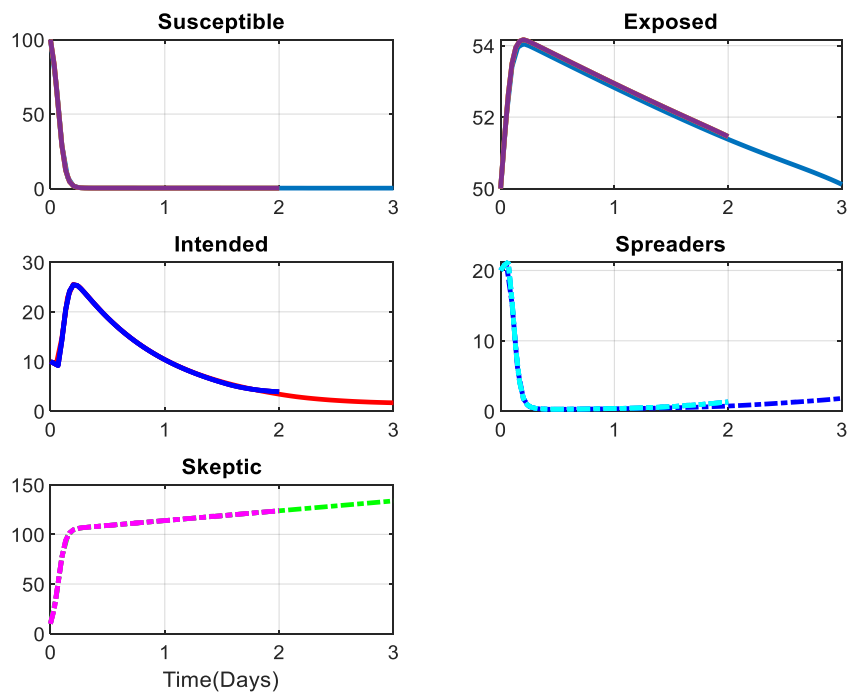

*Figure 6: State trajectories showing how populations changing after 2 and 3 days for the application of control variables*

Figure 6 shows a sharp decrease of Exposed from approximately 54 to around 50 within these days, minimization of Intended/Source by 1 days and almost elimination by 3 days, Mediums shows a constant behavior till 2nd days and a sharp increase after the day 2, Skeptics/Protesters shows a linear increment over time.

The rapid increase in Skeptics meaning the effectiveness of the disinformation campaign in generating skepticism. This aligns with the rate-induced tipping (R-tipping) concept, where the rate of change, not just the total amount of disinformation, can trigger tipping points.

The far-reaching implications of this study extend beyond theoretical understanding. Our findings pave the way for developing a new generation of regulation measures and early warning systems. By monitoring network dynamics and the rate of change in external forces, we can identify nascent disinformation campaigns before they gain traction. This paves the way for targeted fact-checking initiatives, media literacy programs, and strategic counter-narratives aimed at specific audience segments within the network. In conclusion, this study offers a powerful lens through which to decode the enigmatic dance of disinformation. By harnessing the principles of R-tipping and network analysis, we can equip ourselves with the tools to predict, control, and ultimately overcome this dangerous threat. By ensuring the flow of accurate and reliable information, we can safeguard democratic processes, foster informed decision-making, and build a more resilient future for all.